# “Joint Optimization of Opportunistic Predictive Maintenance and Multi-location Spare Part Inventories for a Deteriorating System Considering Imperfect Actions”

## Abstract

Considering the close interaction between spare parts logistics and maintenance planning, this paper presents a model for joint optimization of multi-location spare parts supply chain and condition-based maintenance under predictive and opportunistic approaches. Simultaneous use of the imperfect maintenance actions and innovative policy on spare part ordering, which is defined based on the deterioration characteristic of the system, is a significant contribution to the research. This paper also proposes the method to determine the inspection time which not only considers restraints of the both maintenance and spare parts provision policies, but also uses an event-driven approach in order to prevent unnecessary inspections. Defined decision variables such reliability, upper limit for spare parts order quantity, preventive maintenance threshold, re-ordering level of degradation, and the maximum level of successive imperfect actions will be optimized via stochastic Monte-Carlo simulation. The optimization follows two objectives: (1) system should reach the expected availability which helps decision makers apply the opportunistic approach (2) and cost rate function as an objective function must be minimized. To illustrate the use of the proposed model, a numerical example and its results finally is presented.

**Key words:** maintenance, multi-location supply chain, spare parts inventory, imperfect maintenance, predictive inspection, opportunistic approach, availability, reliability

**Morteza Soltani**

Graduate student of Industrial Engineering,

Iran University of Science and Technology,

Email: m_soltani@ind.iust.ac.ir

Webpage: www.morteza-soltani.ir

## Nomenclatures

| | |
|---|---|
| $X_i$ | system deterioration level at i[th] inspection |
| $Qu_i$ | ordering quantity at i[th] inspection |
| $Tr_i$ | delivery time of the order at i[th] inspection |
| $\alpha_k$ | scale parameter of deterioration process after the k[th] imperfect maintenance action |
| β | shape parameters of the deterioration process when the system is as good as new |
| $\nu_k$ | mean deterioration speed after the k[th] maintenance action |
| $L$ | failure threshold |
| $Z^k$ | k[th] intervention gain |
| $\gamma$ | non-negative real number and represents the impact of imperfect maintenance actions on the deterioration speed of the system |
| $\eta$ | a non-negative real number |
| $T_i$ | i[th] inspection time |
| $\varepsilon_k$ | speed acceleration of deterioration process after k[th] imperfect action |
| $N_{ins}(t)$ | number of inspection in [0,t] |
| $N_{ip}(t)$ | number of imperfect actions in [0,t] |
| $N_c(t)$ | number of corrective actions in [0,t] |
| $N_o(t)$ | number of ordinary orders in [0,t] |
| $N_{oe}(t)$ | number of emergency orders in [0,t] |
| $LT_{s1}$ | lead time of local supplier 1 |
| $LT_{s2}$ | lead time of local supplier 2 |
| $LT_{se}$ | lead time of main supplier |
| $P_{s1}$ | probability of spare part provision by local supplier 1 |
| $P_{s2}$ | probability of spare part provision by local supplier 2 |
| $P_{se}$ | probability of spare part provision by main supplier |
| $C_{ins}$ | cost of each inspection |
| $C_h$ | holding cost rate |
| $C_{oe}$ | emergency ordering cost |
| $C_c$ | corrective maintenance cost |
| $C_{d_1}$ | malfunction cost rate |
| $C_o$ | ordinary ordering cost |
| $C_{pur}$ | purchasing cost of each spare part |
| $C_{d_2}$ | downtime cost rate |
| $d_1(t)$ | malfunction time |
| $d_2(t)$ | downtime cost |
| $C_P^k$ | cost of the k[th] imperfect maintenance action |
| $C_P^0$ | perfect preventive maintenance cost |
| $t$ | simulation time |
| $CMS$ | needed spare part for corrective maintenance |
| $PMS$ | needed spare part for preventive perfect maintenance |
| $IPMS$ | needed spare part for preventive imperfect maintenance |
| $M$ | preventive maintenance threshold |
| $K$ | imperfect maintenance threshold |
| $T$ | re-order level of degradation |
| $S$ | maximum level of spare part inventory |
| $Q$ | failure probability between two inspection times |

# 1- Introduction

Deteriorating systems are known as the systems that deteriorate gradually and the degradation continues until the system ends in a failure. The degradation process usually is presented by failure threshold which directly depends on the maintenance approach used for the system. There are many approaches applied for deteriorating systems' maintenance especially preventive approaches; such as periodic, age-based and condition-based maintenance. Thanks to the rapid development of monitoring equipment which provides accurate information about the system condition, CBM becomes nowadays more and more popular approach in industrial application.

Although maintenance approach is the first factor which crosses the mind at first glance, spare parts provision and inventory policies are the other crucial factors that cannot be ignored in operations optimization. The periodic and continuous are the two frequent kinds of policies considered for spare parts ordering. There are interesting concepts which have significant effects on joint optimization of maintenance and spare parts provision like imperfect maintenance, failure modes, delay time and spare parts deterioration.

Spare parts supply chain is the other aspect of joint models that has direct influence on the shortage or surplus of the stocks. Lead time, which is under the effect of supply chain structure such as multi location, multi echelon, multi indenture etc., is an important factor that can be considered as either deterministic or stochastic parameter.

To find out the optimal values of variables and parameters, indexes including reliability, availability and cost rate are frequently used as the main criteria. Contemplating all aforementioned aspects, joint optimization of maintenance planning and spare parts provision closely depends on an integrated approach which enjoys both inventory and maintenance policies at the same time in order to develop a much more effective decision support system used by beneficiaries and final decision makers. The scope of the joint model has been shown in figure 1.

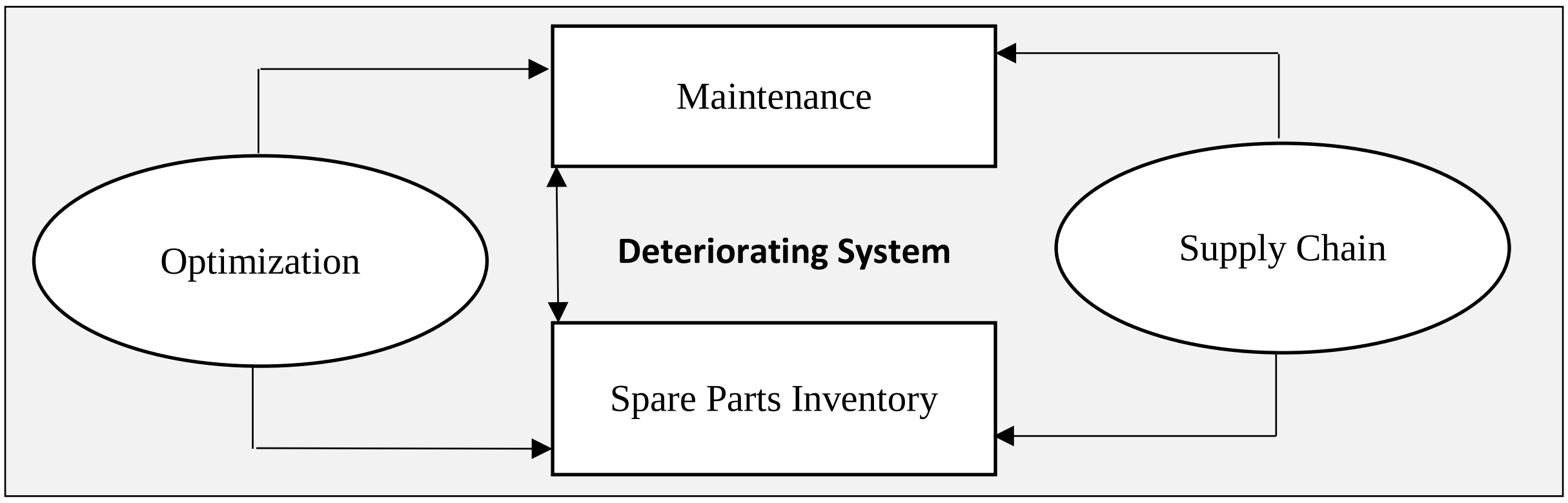

*Figure 1 – Scope of the problem*

The rest of the paper is organized as follows. Section 2 reviews the literature of joint models for deteriorating systems. Section 3 is devoted to the descriptions of the characteristics of the system that should be optimized and other related assumptions. An adaptive joint model is developed and described in Section 4. To evaluate the proposed model, some numerical values of the system's parameters are introduced and related numerical results are in addition discussed in Section 5. Finally, the last section presents the conclusions drawn from this work.

# 2- Literature review

In this section we will review the literature of "Joint Optimization Models for Deteriorating Systems" which mainly consists of three fundamental research areas;

1- Systems degradation
2- Maintenance approaches
3- Inventory policies

Pierskalla and Voelker [37] conducted a survey on maintenance models for deteriorating systems and McCall applies maintenance policies for stochastically failing equipment [6,28]. Subsequently, other researchers use new maintenance concepts for optimization of the deteriorating systems.

Falkner [16] was the first one who presented "joint optimization" title by combining both maintenance and inventory policies. His method for developing the joint models was dynamic programming which could help him find the initial solution for spare part inventory and define the interval time of inspections. Osaki et al. [36] focus on ordering policy with lead time to find the optimal time of ordering and replacement by optimization of the predefined inventory functions. Researches about joint models could be classified in 6 main classes based on the related maintenance and spare parts inventory criteria. This classification has been shown in Table 1. [43]

***Table 1 – Joint Models Classification***

| Class NO. | First criteria | Second criteria | Percent of published papers for each class |
|---|---|---|---|
| | Spare parts inventory policy | Maintenance policy | |
| 1 | Periodic | Periodic | 20 % |
| 2 | | Age-based | 6 % |
| 3 | | Condition-based | 12 % |
| 4 | Continuous | Periodic | 16 % |
| 5 | | Age-based | 10 % |
| 6 | | Condition-based | 36 % |

Comparing the results of joint and sequential models, Acharya et al. [2] proved that joint models are more effective and reliable which resulted in much more extension and development of joint models. Acharya et al. and Chelbi and Aït-Kadi, applied periodic approaches for both maintenance and inventory policies, but the Chelbi and Aït-Kadi used numerical calculation instead of simulation.[8]

In this regard, the two papers published by Amstrong and Atkins [3,4], that are placed in second class in Table 1, consider a single unit system. However, considering operational costs, service restraints and random lead time makes one of them more complex than the other. As the complexity of this type of models has been increased, using simulation-based models have been recently developed intensively. Kabir and Al-Olayan's articles [1,53], both of which are classified in fifth class, conducted simulation-based optimization for both single-unit and multi-units systems. Application of genetic algorithms in simulation process should be considered as an important contribution to joint models, as Ilgin and Tunali [20] use these algorithms to solve the complex descriptive models. Considering downtime costs caused by spare parts shortage, Sarker and Haque [39] conducted a survey on joint models applying for production systems. Nguyen and Bagajewicz [31] developed maintenance model for a processing plants by using the concepts of the failure modes and human resources constraints, and investigated on their effects on system's malfunction or downtime.

During the time that joint models have been developing, several optimization methods were used in deteriorating systems. On the other hand, technology development provided an opportunity to monitor systems' conditions and led decision makers to the vast use of CBM for deteriorating systems. [42]

Using CBM for a single-unit deteriorating system, Grall [18] employed preventive maintenance threshold for predicting the inspection time. Elwany and Gabraeel [14], emphasizing on the effects of technology on maintenance and monitoring, not only developed an integrated model for condition-based (sensors-driven) maintenance and spare part inventory but also employed the concept of RUL in predictive maintenance. Therefore, Elwany's paper could be considered as a connecting point between joint models and deteriorating systems. Most of Wang's works focus on the optimization of maintenance and spare parts inventory for deteriorating systems which is closely related to the topic of this research. [45] Considering a deterioration level for spare parts reordering is one of his innovative idea about inventory policy [26]. Moreover, Wang employs Markov process for optimizing a multi-components system and uses reliability and availability as two main optimization objectives in order to find the optimal inspection time. [46]

There are various rules and policies which are specific to each system with specific characteristics and designs. Thus, the developed theories and models in this area of research are too complicated to solve in an exact way. 'Optimization via simulation' is the frequent method for solving this kind of problems. Figure 2 illustrates the development trend of the joint models used for deteriorating systems. To review and compare all aspects of the related papers, the review table is illustrated in Table 2.

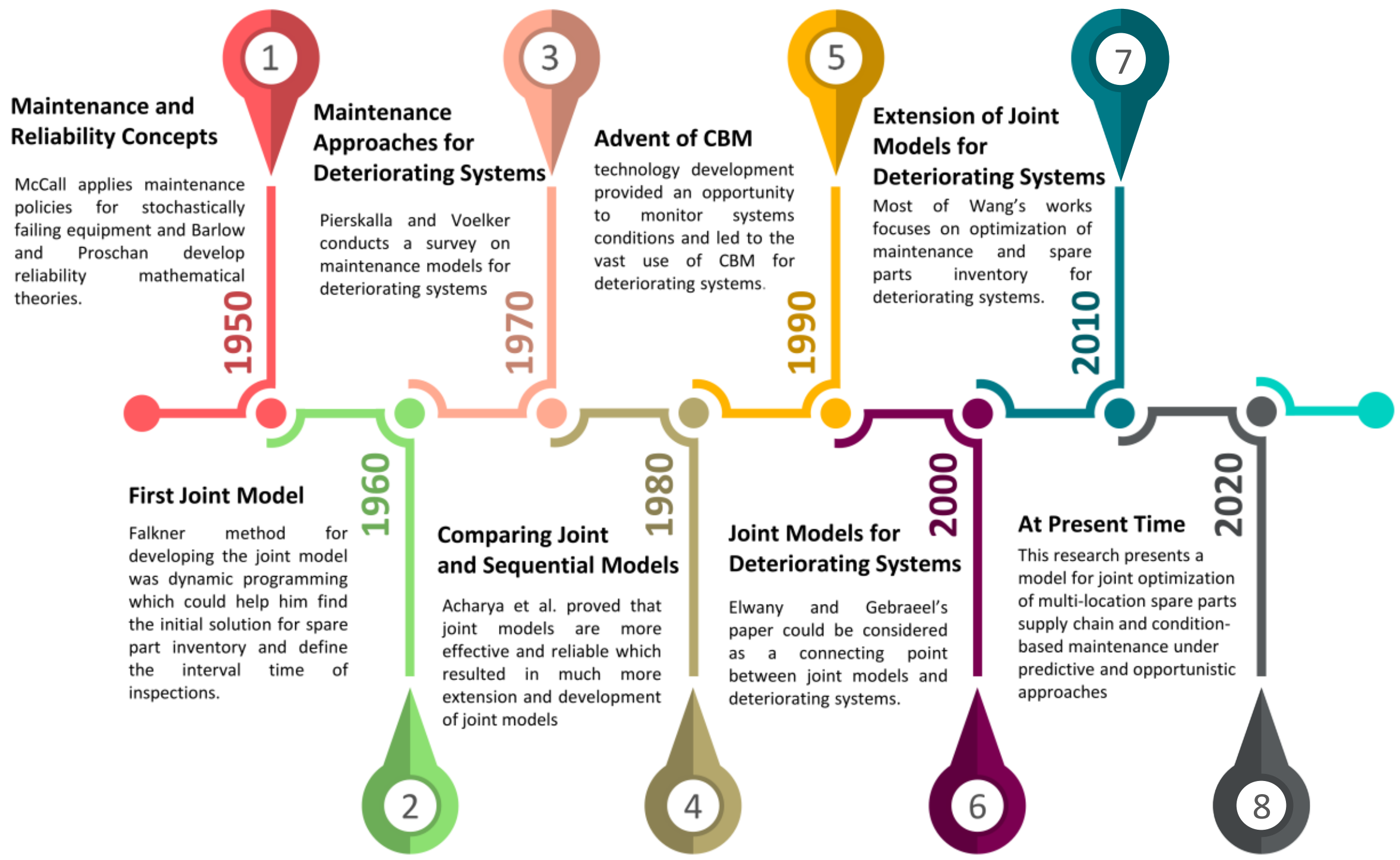

*Figure 2 – the development trend of the joint models used for deteriorating systems*

*Table 2 – Review Table*

| NO. | Co-Author | Year | Deteriorating System? | Joint Model? | Maintenance | | | | | Inventory | | | Optimization | | | Supply chain | | |
|---|---|---|---|---|---|---|---|---|---|---|---|---|---|---|---|---|---|---|
| | | | | | Frequent Policy | | | Innovation | | Frequent Policy | | Innovation | Objective function | Innovation | | Lead time | | Innovation |
| | | | | | Periodic | Age-based | Condition based | Opportunistic | Imperfect maintenance | Periodic | Continuous | Degradation Reorder level | Cost rate | Availability | Reliability | Random | Constant | Multi location |
| 1 | Pierskalla | 1976 | ✓ | × | ✓ | ✓ | | | | | | | | | | | | |
| 2 | Acharya | 1986 | × | ✓ | ✓ | | | | | ✓ | | | ✓ | | | | | |
| 3 | Flores | 1989 | ✓ | × | ✓ | ✓ | | | ✓ | | | | | | | | | |
| 4 | Kabir | 1994 | × | ✓ | | ✓ | | | | | ✓ | | ✓ | | | ✓ | | |
| 5 | Armstrong | 1996 | × | ✓ | | ✓ | | | | ✓ | | | ✓ | | | | ✓ | |
| 6 | Kabir | 1996 | × | ✓ | | ✓ | | | | | ✓ | | ✓ | | | ✓ | | |
| 7 | Armstrong | 1998 | × | ✓ | | ✓ | | | | ✓ | | | ✓ | | | ✓ | | |
| 8 | Sarker | 2000 | × | ✓ | ✓ | | | | | | ✓ | | ✓ | | | ✓ | | |
| 9 | Chelbi | 2001 | × | ✓ | ✓ | | | | | ✓ | | | ✓ | | | | ✓ | |
| 10 | Grall | 2002 | ✓ | × | | | ✓ | | | | | | ✓ | | | | | |
| 11 | Ilgin | 2006 | × | ✓ | ✓ | | | | | | ✓ | | ✓ | | | ✓ | | |
| 12 | Elwany | 2008 | ✓ | ✓ | | | ✓ | | | ✓ | | | ✓ | | | | ✓ | |
| 13 | L.Wang | 2008 | ✓ | ✓ | | | ✓ | | | | ✓ | ✓ | ✓ | ✓ | | | ✓ | |
| 14 | Nguyan | 2008 | × | ✓ | ✓ | | | | | | ✓ | | ✓ | | | | ✓ | |
| 15 | L.Wang | 2008 | ✓ | ✓ | | | ✓ | | | | ✓ | | ✓ | ✓ | | | ✓ | |
| 16 | L.Wang | 2009 | ✓ | ✓ | | | ✓ | | | | ✓ | | ✓ | ✓ | ✓ | | ✓ | |
| 17 | W.Wang [47] | 2011 | × | ✓ | ✓ | | | | | ✓ | | | ✓ | | | | ✓ | |
| 18 | Costantino [10] | 2013 | × | ✓ | | | ✓ | | | | ✓ | | ✓ | ✓ | | ✓ | | |
| 19 | Zanjani [23] | 2014 | × | ✓ | ✓ | | | | | | ✓ | | ✓ | ✓ | | | ✓ | |
| 20 | Wang [48] | 2015 | ✓ | ✓ | | | ✓ | | | | ✓ | ✓ | ✓ | | | ✓ | | |
| 21 | Jiang [21] | 2015 | × | ✓ | ✓ | | | | | ✓ | | | ✓ | | | | ✓ | |

| NO. | Co-Author | Year | Deteriorating System? | Joint Model? | Maintenance | | | | | Inventory | | | Optimization | | | Supply chain | | |
|---|---|---|---|---|---|---|---|---|---|---|---|---|---|---|---|---|---|---|
| | | | | | Frequent Policy | | | Innovation | | Frequent Policy | | Innovation | Objective function | Innovation | | Lead time | | Innovation |
| | | | | | Periodic | Age-based | Condition based | Opportunistic | Imperfect maintenance | Periodic | Continuous | Degradation Reorder level | Cost rate | Availability | Reliability | Random | Constant | Multi location |
| 22 | Phuc Do [12] | 2015 | ✓ | × | | | ✓ | | ✓ | | | | ✓ | | | | | |
| 23 | Kader [22] | 2016 | × | ✓ | ✓ | | | | | ✓ | | | ✓ | | | | | |
| 24 | Shi [40] | 2016 | × | ✓ | ✓ | ✓ | | | | | ✓ | | ✓ | ✓ | ✓ | | ✓ | |
| 25 | Keizer [35] | 2017 | × | ✓ | | | ✓ | | | | ✓ | | ✓ | ✓ | | | ✓ | |
| 26 | Zhang [51] | 2017 | ✓ | × | | | ✓ | ✓ | ✓ | | | | ✓ | | ✓ | | | |
| 27 | Zhang [52] | 2017 | ✓ | ✓ | | | ✓ | ✓ | | ✓ | ✓ | | ✓ | | | | ✓ | |
| 28 | Nguyan [32] | 2017 | ✓ | ✓ | | | ✓ | | | ✓ | | ✓ | ✓ | | | | ✓ | |
| 29 | Zahedi [50] | 2017 | ✓ | ✓ | | | ✓ | | | ✓ | ✓ | | ✓ | | | | ✓ | |
| 30 | Eruguz [15] | 2017 | ✓ | ✓ | | | ✓ | | | | ✓ | | ✓ | | | | | |
| 31 | Siddique [41] | 2018 | ✓ | ✓ | ✓ | | | | | ✓ | | | ✓ | | | | ✓ | |
| 32 | Soltani | 2018 | ✓ | ✓ | | | ✓ | ✓ | ✓ | | ✓ | ✓ | ✓ | ✓ | ✓ | ✓ | | ✓ |

# 3- System description and assumptions

In this section we will focus on the system's characteristics, necessary assumptions and rules which should be considered to define the descriptive problem for a specific type of deteriorating asset. In this system, there is no way to monitor the sub-system's failure signs except operational inspection. Thus, inspection would be the only way to get information about the systems' degradation level.

According to the random nature of the failure events, the failures could happen between inspection times even if we define a short time interval for periodic inspections. There are many maintenance actions, which could be conducted during the inspection time, three of which I will explore below.

1- **Corrective maintenance** is done after initiation of failure, leading to degraded performance which could be revealed by inspections. In this situation, the system may be led into malfunction or breakdown that will result in high costs. After performing a corrective maintenance action, the system goes back to the 'as good as new' state.
2- **Perfect preventive maintenance** is done before initiation of failure. Alike other types of preventive maintenance, perfect maintenance extends equipment lifetime, but it causes high cost because the deterioration level must be reset to zero. Compared to the corrective maintenance, the preventive perfect maintenance cost is much lower.
3- **Imperfect preventive maintenance** is done before initiation of failure. The system condition after this action will be somewhere between the condition before maintenance and as good as new. From a practical point of view, imperfect maintenance can describe a large kinds of realistic maintenance actions. The imperfection of these actions arises from two main reasons. The first one is the "bad" realization of a perfect maintenance action due to human factors (e.g. stress, lack of skills, lack of attention), shortage of spare parts, lack of repair time, etc. The other reason could be considered as a deliberate one, for instance, the maintenance policy of decreasing costs which may lead to deal with "low-cost" people, spare parts, and logistics. While the first reason does not provide any benefit, the second one may lead to cost benefits. On the other hand, accumulations of short deterioration induce long-term deterioration that could be technically or economically no-more acceptable. [12]

Preserving ample sizes of spare part inventories for immediate disposition whenever needed, can be a logical solution to the spare parts' availability problem. However, this solution entails a high stocking cost. Thus, there must be a trade-off between overstock and shortages of spare parts which is an inventory planning problem with a maintenance scheduling aspects. Spare parts procurement depends directly on the structure of supply chain. In this problem, there are two local suppliers and one main supplier each of which has its own specific lead time that has an effect on ordering priority. So, the supplier with less lead time would be prior to order. Moreover, there is a possibility to order from the main supplier in an emergency condition if none of the local suppliers would not be able to supply the needed spare parts.

As the system's reliability is adaptive in this problem, by modifying the quality of parts or skill of technicians, decision makers are able to adapt the level of reliability to reach the optimal state of the system. There are some significant questions that this research aims to find their answers;

- How to determine the next inspection time?
- What is the best time to order and what is the optimal quantity of the order?
- How to choose the best supplier?
- What is the optimal reliability of the system?
- What maintenance approaches should be used and when should we perform each of them?
- To find the optimal value of decision variables, what are key indexes for optimization?
- How to solve this complex problem?

There are several assumptions that have been made in this problem;

- We assume that all imperfect actions will be performed because of the second reason in order to provide costs benefits.
- Imperfect maintenance has two effects on the deteriorating system. First, imperfect maintenance restores a system to a state between good-as-new and bad-as-old. Second, each imperfect preventive action will accelerate the speed of the system's deterioration process.
- Imperfect maintenance carries less costs in comparison with corrective and perfect maintenance. The perfect maintenance needs more accuracy to be performed. The more accuracy is required for the action, the more time and cost will be needed. The other reason for the higher cost of perfect maintenance is using much more resources in compare to imperfect one.
- The intervention gain in the deterioration level of the system due to an imperfect maintenance action is random. The cost of each imperfect action is in proportion with the deterioration improvement. From a practical point of view, in most cases, the quality of the maintenance action increases with the level of resources allocated to it, and hence with its cost. [9,27,30]
- During the imperfect maintenance, spare part replacement may happen which means that replacement would not necessarily lead the system to the state of 'as good as new'. The important point to be mentioned is that the quality of maintenance depends on not only the spare parts replacement but also the quality of performing the related maintenance tasks.
- The lead time and transport cost depend on the distance of operation place from suppliers. Therefore, lead time could be considered constant and correlated to the distance between suppliers and operations site.
- We assume that there is no price difference between the spare parts supplied by different suppliers.
- Inspection is the only way to be aware of the malfunction of the system.
- Considering the failure threshold, the system will start malfunctioning if the deterioration level of the system passes over the failure threshold.
- The deteriorating system is the critical sub-system whose downtime will definitely lead production system to downtime.
- The system is simplified by a single component system (e.g. only the most important component is considered)
- It is assumed that maintenance duration are negligible.
- We will consider loss benefits caused by system's downtime or malfunction as the opportunity costs.
- Emergency ordering cost (from main supplier) and ordinary ordering cost (from local supplier) are different. Moreover, the ordering costs do not depend on the quantity of each order because these costs closely relate to the procurement process including employees' salaries.
- We consider the constant overhaul interval time for the production system which is determined by decision maker based on overhaul policy for the total system. The system's scheduled downtime (for overhaul) will provide an opportunity to shift the adoptive corrective maintenance time to the prescheduled overhaul time. Simultaneous performing corrective maintenance of all sub-systems will lead the system to a considerable cost saving because it prevents unnecessary total system's downtime.

# 4- Model development

In this section, we will develop an integrated model based on the related assumptions to optimize the maintenance planning of the deteriorating system and spare part provisioning jointly. The model consists of 6 main areas which have close interaction with each other (See figure 3).

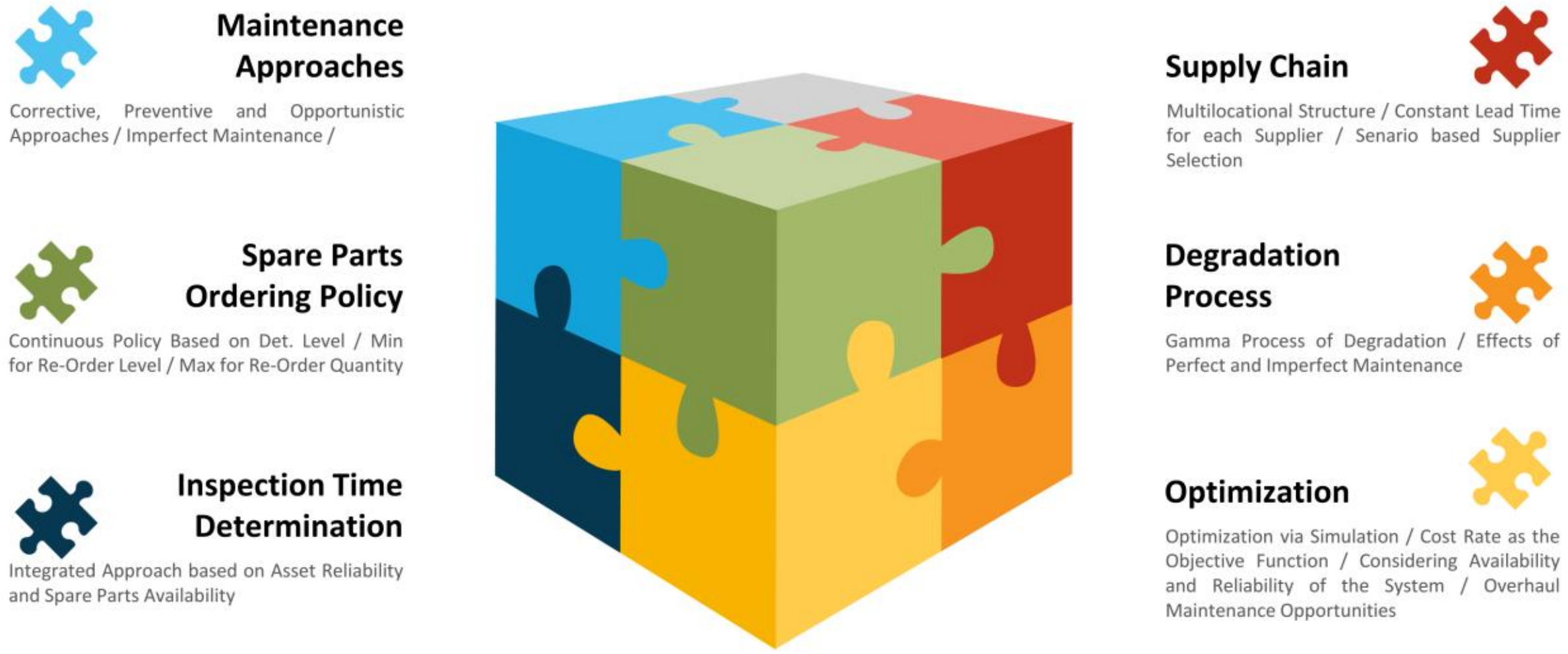

*Figure 3 – 6 main areas of the problem*

## 4-1- Spare parts ordering policy

Despite of the frequent policies used for inventory optimization, which consider reorder point, we employ reorder level based on deterioration characteristic of the system. Thus, we defined the lower level of degradation (T) [45]. At each inspection time, we will order if the degradation level after performing the expected operation is upper than the predefined T level. However, the important question to be answered is how the optimal quantity of the order is calculated. At each decision point, we order the spare part up to the upper limit of order quantity (S) which will be optimized through simulation as one of our decision variables. In this model, total stock of the spare part is equal to the sum of stock in hand and undelivered stock ordered previously. At each decision point:

Qu = S – Total stock (1)
Where: Total stock = Stock in hand + Undelivered stock (2)

As a result, we define (T,S) policy for the spare part inventory which means we will order up to S if the degradation level of the system after probable maintenance action is upper than T level. In this model, ordering policy and inspection time will be defined so that the shortage would not happen at the inspection times before the system's failure. The reason is that we employ an integrated approach which considers the time of having enough spare part for expected maintenance operations and the time of system's failure at the same time. Thus, at each decision point, we should determine both in hand and undelivered stocks for the next inspection time. We define Tr to predict the delivery time of the order based on the lead time coordinated with selected supplier.
On the other hand, we must order instantly if we face shortage at failure time. In other words, if the next inspection time is in the failure zone (degradation level passes over failure threshold), we should place an order instantly. The order quantity must be defined so that the stock quantity would be at maximum level of S after performing the corrective maintenance. At this situation, the next inspection time is the time that provides enough spare parts for conducting the corrective maintenance.

### 4-2- Degradation process

Defined as the duration left for a system before it fails, Residual Useful Life (RUL) has been recently introduced especially for predicting the inspection time [11,17,49]. Gamma processes have been widely used to describe the degradation of systems [19,29,44]. Strictly monotone increasing, which is the behavior observed in most physical deterioration processes, is the characteristic that justifies using Gamma process. It is assumed that the system's deterioration between the k[th] and the (k+1)[th] maintenance actions evolves like a Gamma stochastic process. Thus, $(X_{t_2} - X_{t_1})$ as a stochastic variable follows a Gamma probability density (pdf) with shape parameter $\alpha_k(t_2 - t_1)$ and scale parameter β.

$$f_{\alpha_k(t_2-t_1),\beta}(x) = \frac{1}{\tau(\alpha_k(t_2-t_1))}\,\beta^{\alpha_k(t_2-t_1)} x^{\alpha_k(t_2-t_1)-1} e^{-\beta x} I_{(x\geq 0)} \quad (3)$$

$I_{(x\geq 0)}$ is an indicator function where:

$$if\ x \geq 0\,, \qquad I_{(x\geq 0)} = 1$$

$$O.W \quad ,\ I_{(x\geq 0)} = 0$$

$$\alpha_k = \frac{\nu_k}{\beta} \quad (4)$$

$\nu_k$ is the mean deterioration speed of the system between the k[th] and the (k+1)[th] imperfect actions. The system degradation behavior and corresponding states are illustrated (see figure 4)

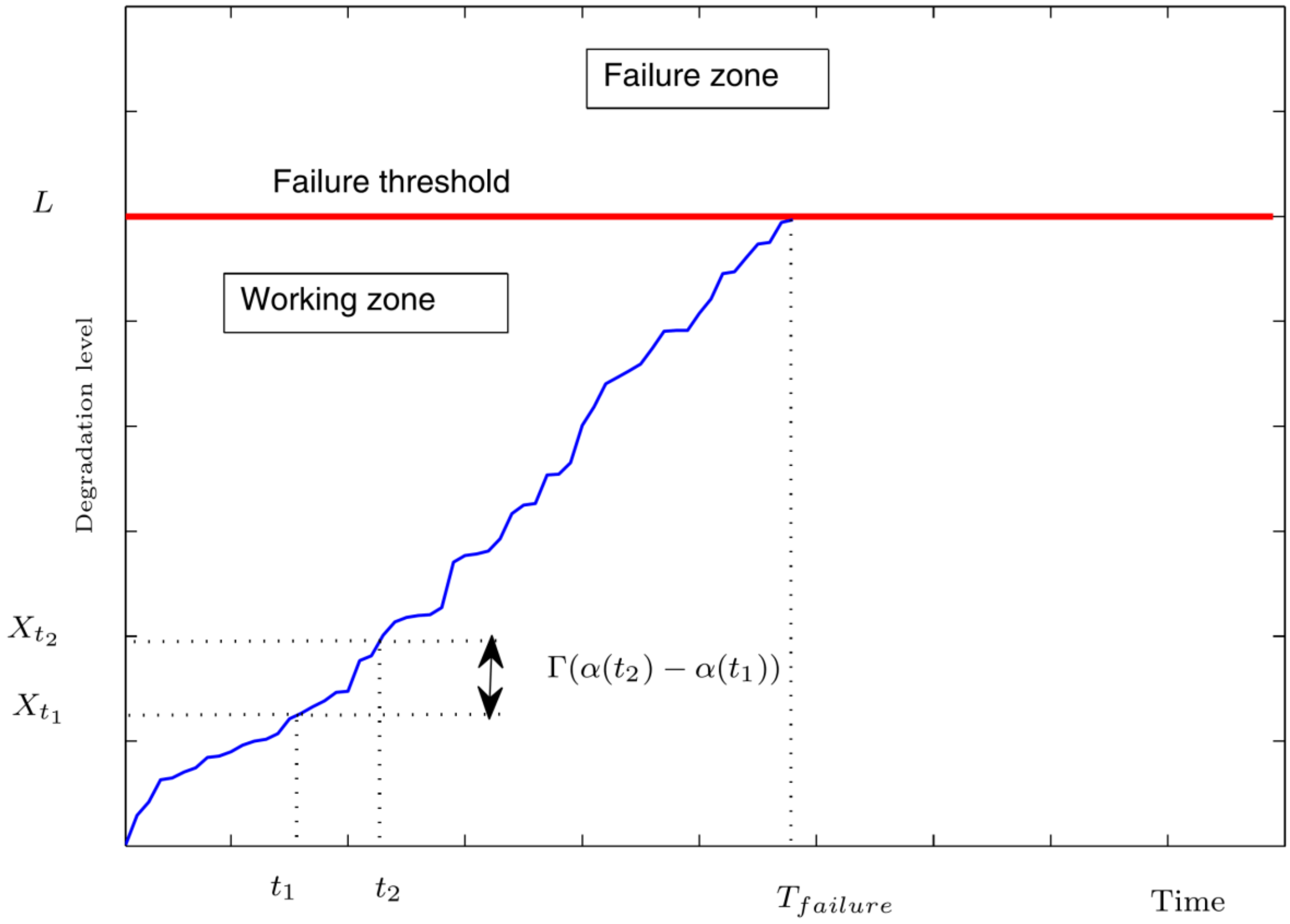

*Figure 4 – system degradation behavior*

## 4-3- Inspection time determination

In this regard, several papers propose event-driven approach to determine the inspection time [5,38], but using policies like periodic inspection time with the fixed interval causes unnecessary inspections carrying a high cost. Therefore, we develop an integrated method to predict the next best time for inspection based on real-time information in order to reduce the unnecessary inspections caused by the periodic approach [13,34]. This prediction is based on all aspects considered in this problem. For this reason, we call this method "Joint approach of inspection determination" which determines the next inspection time at each decision point considering all below questions simultaneously.

1- What is the degradation level of system at this point?
2- How many successive imperfect maintenances have been done?
3- What is the expected reliability of the system?
4- What are the next expected operations?
5- Are there enough spare parts to perform the expected operation?

The main idea of the RUL based inspection is that the next inspection time should be chosen such that the probability of the failure of the system before the next inspection remains lower than a limit Q ($0 \le Q \le 1$). Q is a decision variable to be optimized [17]. According to this inspection policy, the reliability of the system between two inspection times interval remains higher or equal to (1-Q). The inspection time determination includes two steps. First, we use the method 'RUL based inspection'. If we let $T_i$ denote the time at which the system is inspected, and the corresponding degradation level of the system is $X_{T_i}$ (it is the deterioration level of the system after maintenance if a maintenance action is executed at time $T_i$), the next inspection time is then determined by:

$$T_{i+1} = T_i + m\left(X_{T_i}, Q\right) \quad (5)$$

$$m\left(X_{T_i}, Q\right) = \left\{\Delta T : P\left(X_{T_i+\Delta T} \ge L \middle| X_{T_i}\right) = Q\right\} \quad (6)$$

Where:

$$P\left(X_{T_i+\Delta T} \ge L \middle| X_{T_i}\right) = P\left(X_{T_i+\Delta T} \ge L - X_{T_i}\right) = \int_{L-X_{T_i}}^{\infty} f_{\alpha_k \Delta T, \beta}(x)dx = 1 - \int_{0}^{L-X_{T_i}} f_{\alpha_k \Delta T, \beta}(x)dx \quad (7)$$

It is clear that $m\left(X_{T_i}, Q\right)$ depends on the current degradation level of the system, the failure threshold L and the parameter Q. As the second step, we must check if there is enough spare part to execute the expected maintenance actions or not. According to the value of k (the successive imperfect actions that have done), the expected actions could be both perfect and imperfect maintenance. The next inspection time will be equal or more than the time determined in first step, whenever the enough spare parts for the expected action could be provided.

## 4-4- Maintenance approaches

In this section we will investigate the impacts of all kinds of maintenance actions used in this model.

**Corrective maintenance:** the system needs a corrective maintenance if $X_{T_i} \ge L$, the time that deterioration level of the system passes over the failure threshold. After a corrective maintenance action, the system will be led to the 'as good as new' state and the problem comes back to the initial condition (the deterioration level of the system is zero and the spare part inventory is full up to S). As the system is reset to the initial state after corrective maintenance, we can consider the time between two corrective actions as the life cycle of the systems. As a result, the corrective cost ($C_c$) occurs only one time during on life cycle of the system**.**

**Preventive perfect maintenance:** Preventive maintenance (PM) consists of maintenance activities performed before equipment breaks down, with the intent of keeping it operating acceptably and reducing the likelihood of breakdown. The main purpose of PM is to extend equipment lifetime, or at least the mean time to the next failure whose repair may be costly. Furthermore, it is expected that effective PM policies can reduce the frequency of service interruptions and the many undesirable consequences of such interruptions. Perfect PM restores the degradation level of the system to zero in order to reduce the negative impacts of imperfect maintenance. It might be considered as a theoretical case where all preventive maintenance actions are imperfect. After too many successive imperfect actions the system fails and no more imperfect action could restore the system to the working condition because of the high rate of deterioration. Thus, we propose a hybrid policy in which both perfect and imperfect maintenance actions are considered (see figure 5) [13,19], and define the upper limit (K) for successive imperfect actions as the imperfect maintenance threshold. This would be another decision variable optimized in this model.

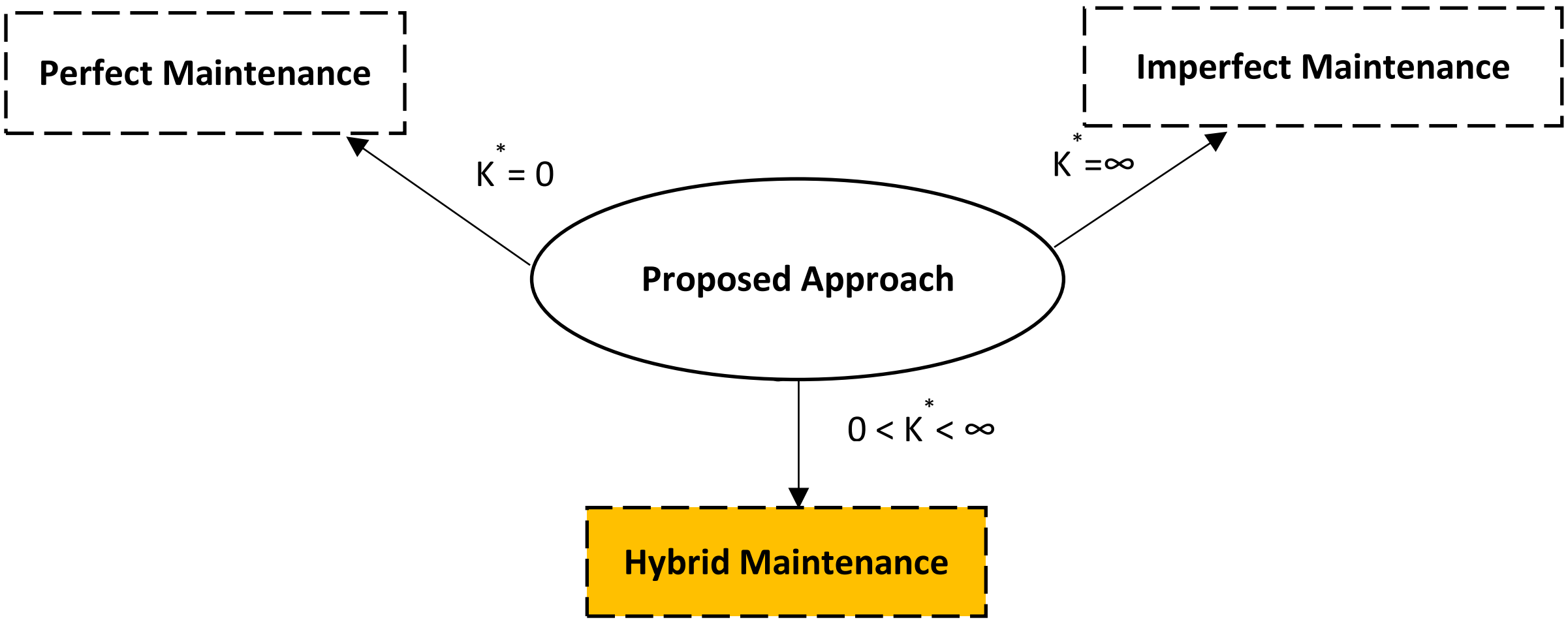

*Figure 5 – Hybrid maintenance approach*

**Imperfect maintenance:** we consider that the imperfect preventive maintenance actions lead the system to a better state for which the degradation level is lower or equal to the current deterioration level of the system. Employing imperfect maintenance makes the model much closer to real conditions although this aspect increases the complexity of the problem [24,33]. The increase in the speed of system's deterioration after each imperfect maintenance can be described by non-negative continuous random variable $\varepsilon_k$ which follows an exponential distribution with density probability:

$$h(x) = \gamma e^{-\gamma x} I_{(x \geq 0)} \; ; \qquad (8)$$

Where:

$$if\ x \geq 0 , \qquad I_{(x \geq 0)} = 1$$

$$O.W , \qquad I_{(x \geq 0)} = 0$$

If the k$^{th}$ maintenance action is a corrective or perfect preventive maintenance, the mean deterioration speed of the system after maintenance is reset to:

$$\nu_k = \nu_0 = \frac{\alpha_0}{\beta} \qquad (9)$$

If the k$^{th}$ maintenance action is an imperfect preventive one, the mean deterioration speed of the system after maintenance is set to:

$$\nu_k = \nu_{k-1} + \varepsilon_k \qquad (10)$$

Figure 6 illustrates the impacts of imperfect maintenance on deterioration process of the system.

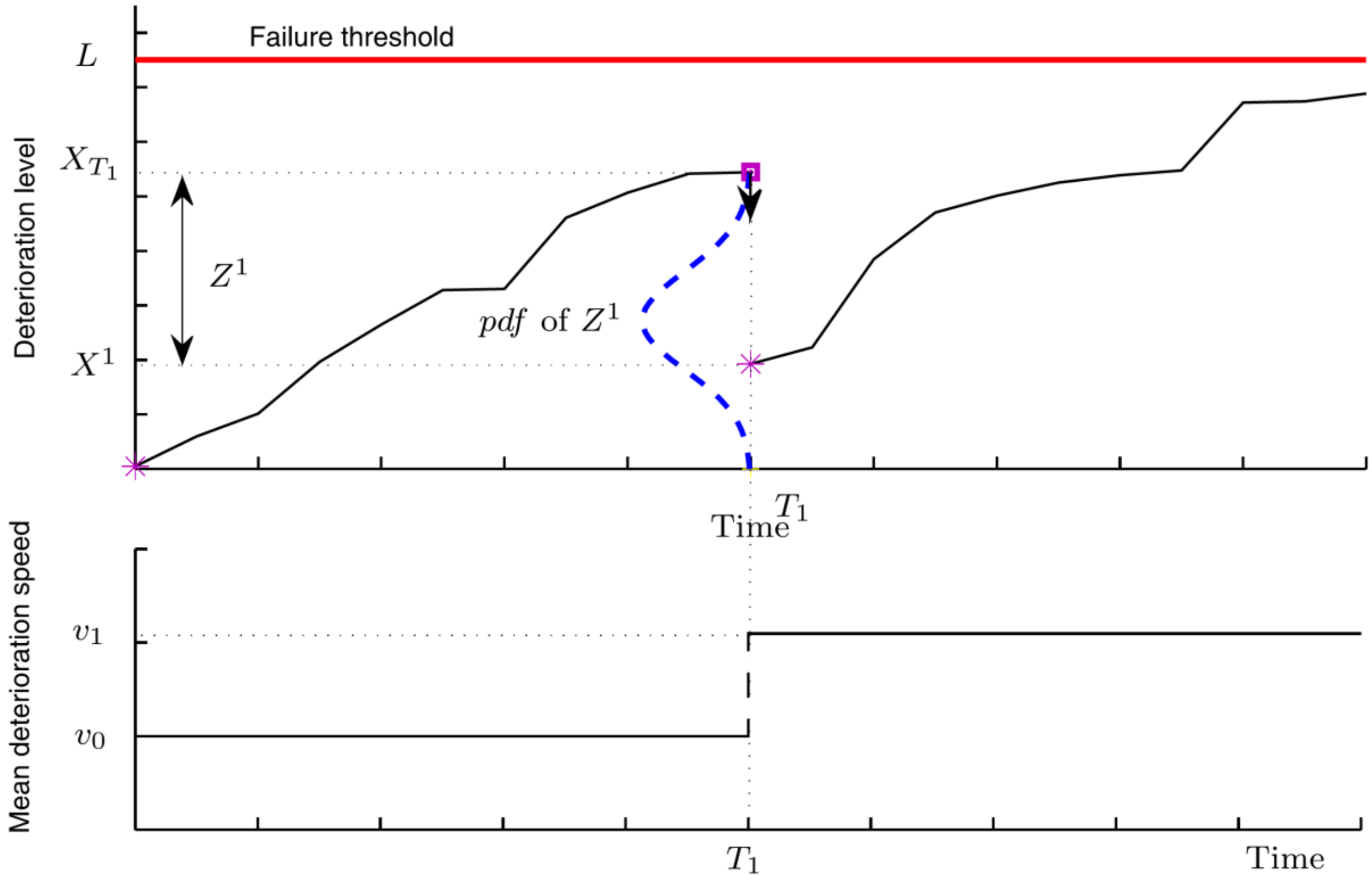

*Figure 6 – deterioration evolution and impacts of imperfect maintenance*

After conducting each imperfect maintenance, the intervention gain is then assumed to be described by a continuous random variable $Z^k$ which distributed with the density [7,29,44]:

$$g_{\mu,\sigma,a,b}(x) = \frac{\frac{1}{\sigma}\phi\left(\frac{x-\mu}{\sigma}\right)}{\Phi\left(\frac{b-\mu}{\sigma}\right) - \Phi\left(\frac{a-\mu}{\sigma}\right)} I_{[a,b]}(x) \quad (11)$$

$$if\ a \leq x \leq b\ , \qquad I_{[a,b]}(x) = 1 \quad (12)$$

$$Otherwise\ \ , \ I_{[a,b]}(x) = 0$$

$\phi(\varepsilon) = \frac{1}{\sqrt{2\pi}\exp(-\frac{1}{2}\varepsilon^2)}$ (13) ; The probability density function of the standard normal distribution

$\Phi(.)$ (14) ; cumulative distribution function

$\mu = \frac{X_{T_i}}{2}$ , $\sigma = \frac{X_{T_i}}{6}$

$a = \mu - 3\sigma$ , $b = \mu + 3\sigma$

For each imperfect intervention, $Z^k$ is bounded, i.e., $0 < Z < X_{Ti}$ where $X_{Ti}$ is the deterioration level of the system at $T_i$.

As we mentioned before, imperfect maintenance carries the least cost among other maintenance actions. The cost of each imperfect action is in proportion with the deterioration improvement. From a practical point of view, in most cases, the quality of the maintenance action increases with the level of resources allocated to it, and hence with its cost. The degradation improvement factor is defined as the ratio of the improvement gain divided by the deterioration level of the system before maintenance [44]. Based on the improvement factor, imperfect maintenance costs can be evaluated and considered as a function of the improvement factor. [9,27,30] Thus;

$$\frac{C_P^k}{C_P^0} = u(T_i)^{\eta} = \left(\frac{Z^k}{X_{T_i}}\right)^{\eta} \qquad (15)$$

It is clear that the perfect maintenance cost ($C_P^0$) is considered as the upper limit of imperfect maintenance cost $C_P^k$ which varies among the boundary ($0 < C_P^k < C_P^0$) depending on the value of $\eta$. Figure 7 illustrates the impacts of $\eta$ on calculation of imperfect maintenance cost and Table 3 elaborates upon the correlation between $\eta$ and $C_P^k$. [7,25]

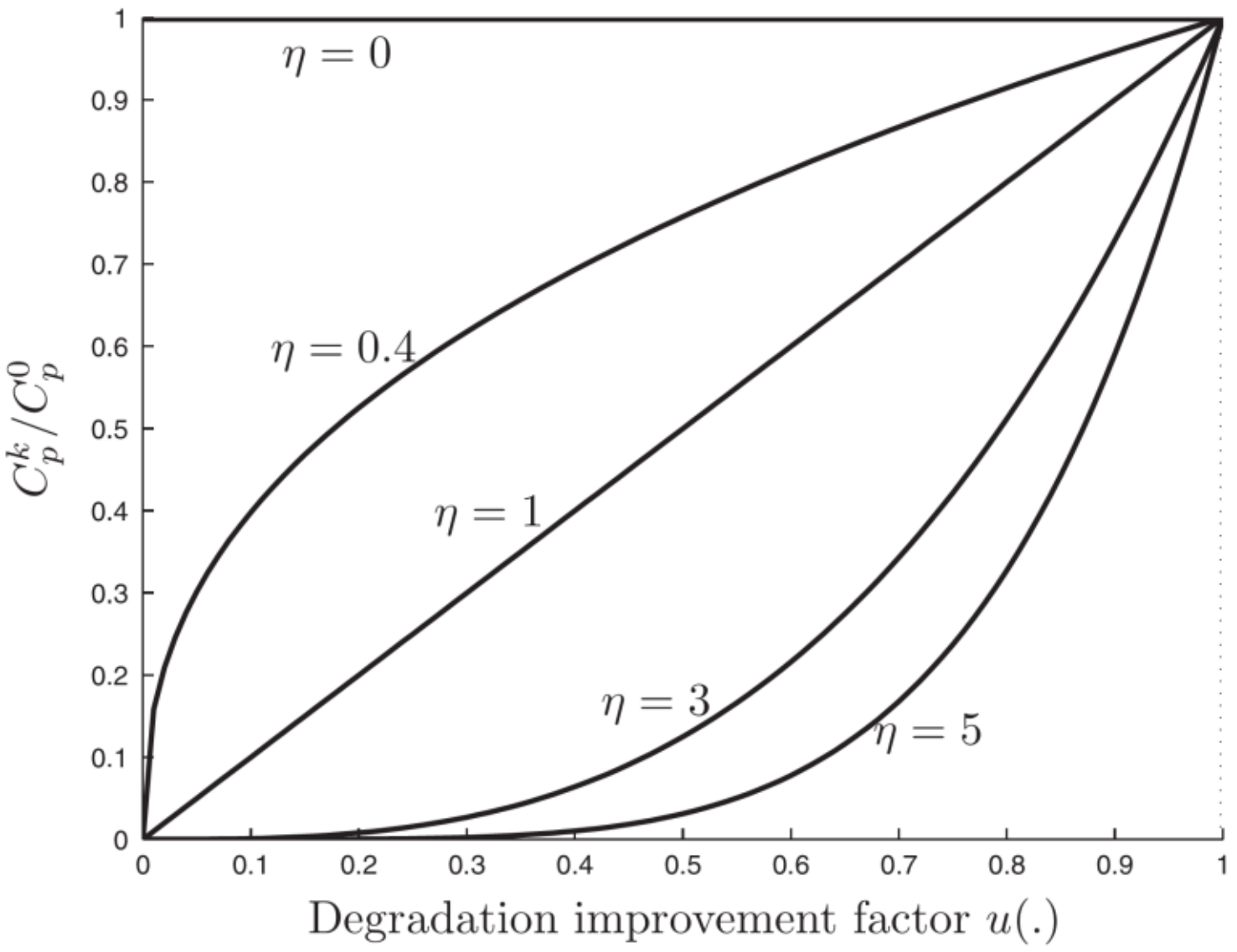

*Figure 7 – imperfect maintenance cost function*

***Table 3 – the correlation between $\eta$ and $C_P^k$***

| NO. | η | $C_P^k$ | Correlation interpretation |
|---|---|---|---|
| 1 | $\eta = 0$ | $C_P^k = C_P^0$ | imperfect maintenance cost is constant |
| 2 | $0 < \eta < 1$ | *$C_P^k$ is a concave function* | the maintenance cost increases more than the improvement gain when performing the maintenance |
| 3 | $\eta = 1$ | *$C_P^k$ is a linear function* | maintenance cost is proportional to the improvement level gain |
| 4 | $\eta > 1$ | *$C_P^k$ is a convex function* | the maintenance cost increases less than the improvement gain |

## 4-5- Supply chain

We assume that there is no price difference between the prices offered by each supplier, so the procurement cost depends on the lead time which is determined based on the distance between the supplier and operation site. It is clear that we prefer to order from the local supplier that is closer to the site whose lead time is less than the lead time of the other supplier. As it is shown in figure 8, there are two local suppliers that one of them provides the ordered spare part sooner, because this supplier is much closer to the operation site. The main supplier supplies both of the local suppliers. Furthermore, there is a possibility to order from the main supplier in an emergency condition if none of the local supplier would not be able to supply the needed spare parts ordered for performing the corrective maintenance. As a result, the emergency order from main supplier occurs when we face shortage in stock for conducting the corrective maintenance, otherwise local suppliers are the only sources that can supply.

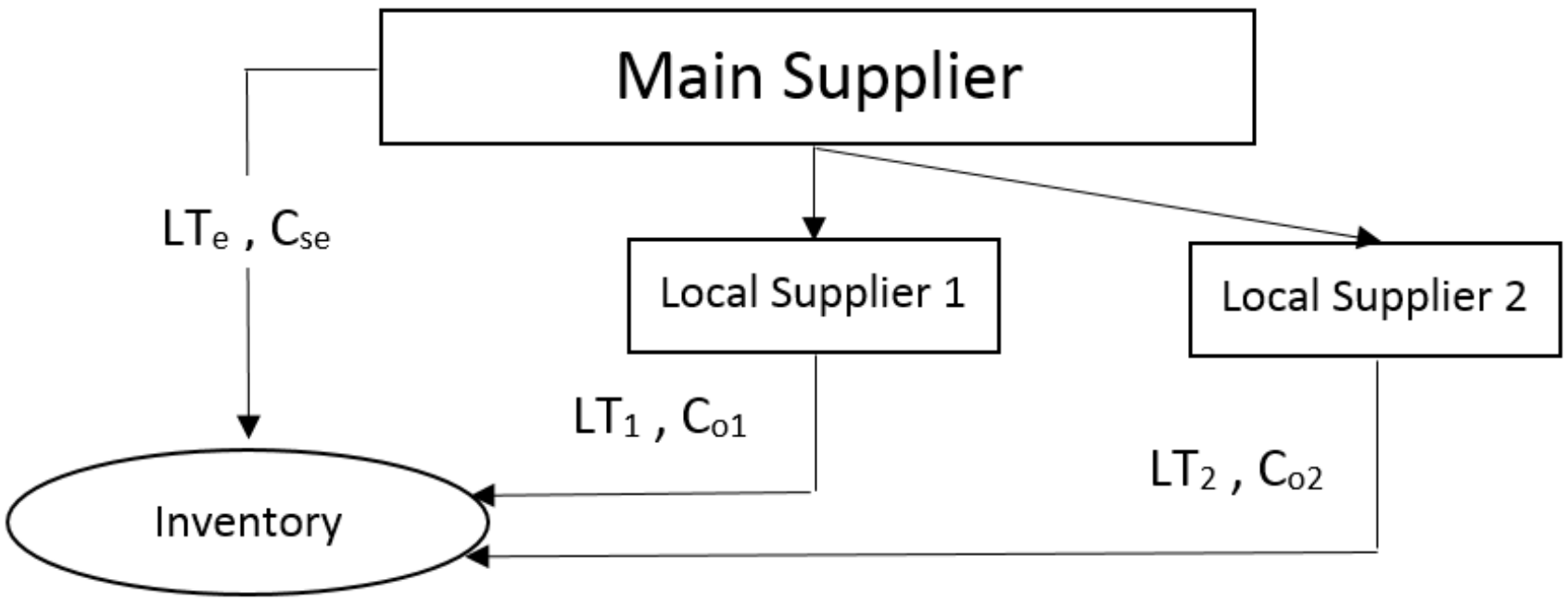

*Figure 8 – supply chain structure*

It is clear that ordering form the supplier with longer lead time causes higher cost of procurement which mainly includes transport costs.

$$LT1 < LT2 << LTe \longrightarrow Cs1 < Cs2 << Cse \quad (16)$$

In practical point of view, we do not know anything about the availability of the spare part at each local supply center. At each decision point, we can only enquire about the availability of needed spare parts. Therefore, the companies are able to calculate the chance of the availability of enough spare part based on the historical data. In order to simulate this condition, we produce a random number between 0 and 1, then compare the produced number with the probability of spare part availability at each supply center. If the random number is less than the probability, it means that there are enough spare parts to order from the supplier. It is obvious that we must start checking the availability with the supplier with shorter lead time. The exquisite point to be mentioned is that the emergency order will be made when the spare part needed for corrective maintenance could not be provided by none of the local suppliers.

## 4-6- Optimization

There are several measures to compare the system's performance under the different conditions such as reliability, availability, and cost rate. The most popular index used commonly in recent research is cost rate. We will develop the mathematical model whose objective function is cost rate. This rate consists of all types of cost considered in the developed model including;

1- **Inspection cost:** For each inspection, we consider the constant cost ($C_{ins}$). The objective function includes all the inspections happened during one life cycle of the system.
2- **Imperfect action cost:** As it is discussed in previous part, the imperfect maintenance cost ($C_p^k$) depends on some factors such as intervention gain and physical characteristics of the system.

In general, each maintenance action incurs a cost and an imperfect maintenance action often incurs a reduced maintenance cost. All imperfect actions conducted during one life cycle of the system will include in objective cost rate.

3- **Corrective maintenance cost:** The corrective cost ($C_c$) occurs only one time during one life cycle time of the system ($N_c(t) = 1$). The duration between the time that system starts to work and the time that the deterioration level of the system passes over the failure threshold which system needs a corrective action, could be considered as the life cycle time of the system. In other words, life cycle time is the duration between two successive states of 'as good as new'.
4- **Preventive perfect maintenance cost:** For each perfect action, we consider a constant cost ($C_P^0$). The objective function includes all the perfect maintenance executed in one life cycle of the system. It is obvious that perfect actions incur higher cost than that of imperfect actions do.
5- **Downtime cost rate:** since we assumed that the operation time is negligible, downtime will occur once the system needs corrective maintenance and we do not have access to enough spare parts. This cost is exactly equal to the income that the production system could earn during the downtime. We consider the constant rate ($C_{d_2}$) for the downtime cost and include it in objective cost rate.
6- **Malfunction cost:** This cost also will happen only one time during the life cycle of the system. The time between system's failure and the next inspection will be known as the malfunction time. The fixed cost rate that we consider for malfunction time ($C_{d_1}$) is lower than the cost rate of downtime. During malfunction time, the production of the system will be continued, but the low quality of the products, which can incur expenses to the system, may not be acceptable.
7- **Holding cost:** we will consider the constant cost rate ($C_h$) for holding each spare part in the inventory. As the time of holding a part increases, the cost incurred will be increases as well. This cost may include energy consumption and employees' salaries.
8- **Ordering cost:** according to the defined supply chain we consider two types of ordering costs. First, the constant cost for the local suppliers that we consider the same value for all the local suppliers. Second, the constant emergency cost which will occur when we need to purchase the spare part from main supplier.
9- **Purchasing cost:** for each spare part purchased in the life cycle of the system, we will consider the constant cost which is as same as the price of the part.

As defined below, we can formulate the cost rate as the objective function of the mathematical model.

$$C^t(M,K,Q,S,T) = C_{ins}.N_{ins}(t) + \sum_{k=1}^{N_{ip}(t)} C_p^k + C_P^0.N_p(t) + C_c.N_c(t) + C_{d_1}.d_1(t) + C_{d_2}.d_2(t) + C_o.N_o(t)$$
$$+ C_{oe}.N_{oe}(t) + C_h.\int_0^t R(x)dx + C_{pur}.\sum_{i=1}^{N_{ins}(t)} Q(i) \quad (17)$$

Furthermore, we will consider the system's availability as the constraint of this model. During the simulation process, the reliability of the system will be optimized as well. In order to apply the opportunistic approach, we define the A* as the lower limit for system's availability based on the interval time between two successive overhauls. Therefore, we should seek for the condition that make the system work continuously (without any failure or downtime) during the time between two overhauls. This condition helps us shift all possible corrective maintenance (for all sub-systems) to the scheduled overhaul time and reduce the number of downtime events. The mathematical model is defined below.

$$Z = MIN\,\{C^t(M,K,Q,S,T)\} \quad (18)$$

$$St: \quad Availibi\,lity \geq A^* \quad (19)$$

# 5- Solving method

In this section we will design the algorithm for simulation process based on the assumptions and the descriptive model discussed in previous sections (see figure 9). The objective function is a response evaluated by the simulation, so we will use simulation optimization to solve this complex problem [20]. This simulation model can be thought of as a mechanism that turns input parameters into output performance measures such as cost rate, system's availability, and reliability. Stochastic Monte Carlo simulation is used to evaluate the performance measures. To do so, for each value of decision variables (M,K,T,S,Q), the corresponding maintenance cost rate as the objective function and system's availability as the opportunistic constraint are calculated. To ensure that the maintenance cost rate is reasonable and valid, a large number of simulation realizations must be done. Moreover, each realization simulates one life cycle of the system. By varying different values of (M,K,T,S,Q), the minimum maintenance cost rates can be identified. However, the optimal values of the decision variables (M,K,T,S,Q) are the corresponding ones to the minimum maintenance cost rates if the availability of the system satisfies the opportunistic constraint.

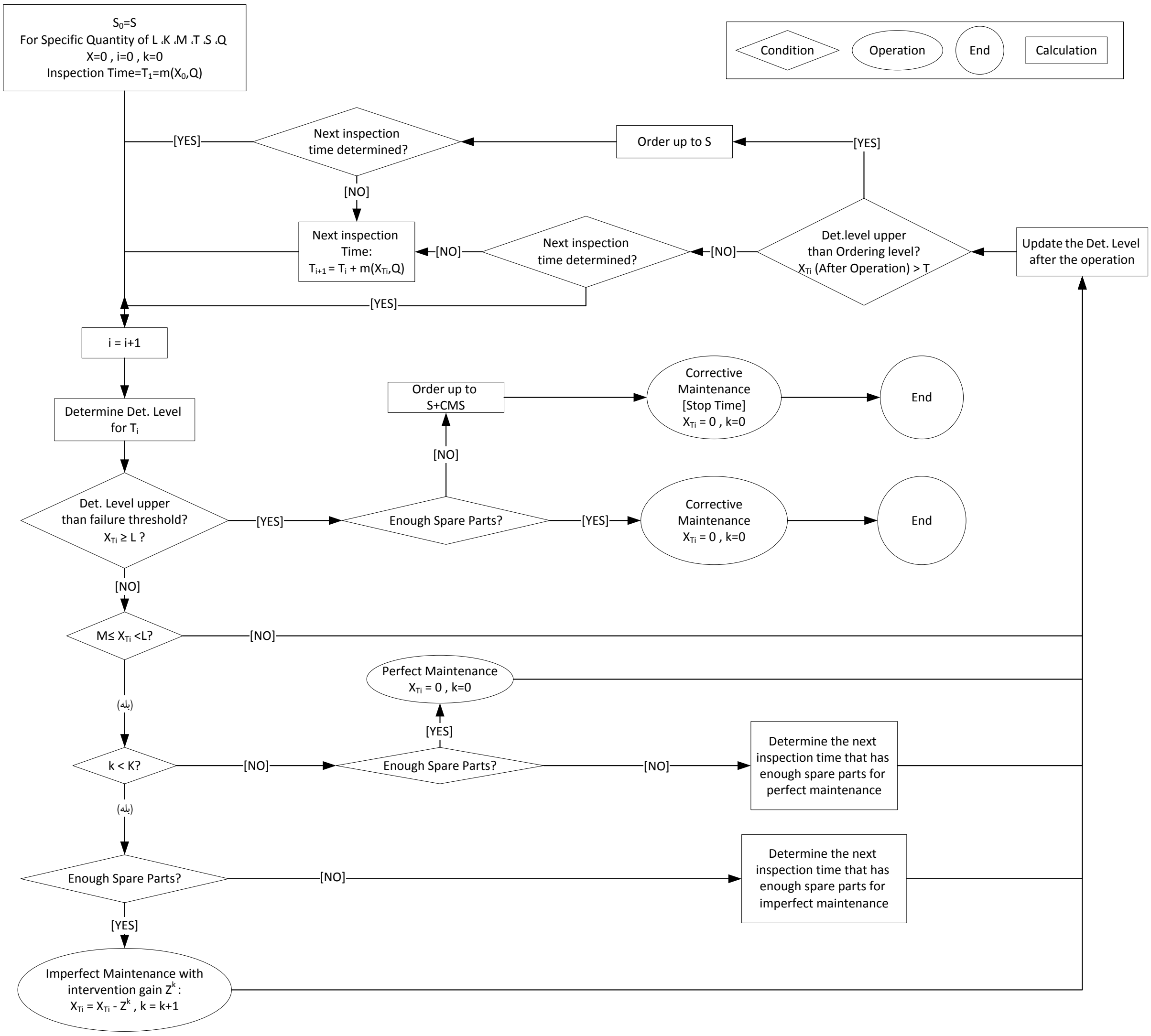

*Figure 9 – simulation algorithm*

## 6- Numerical example and model evaluation

This section intentionally left blank

# 7- Conclusion and future research

In this research, we focus on designing a decision support system for optimizing the maintenance planning and spare parts provision for a deteriorating system. We also apply Monte-Carlo simulation to solve this complex problem. Contemplating all related aspects, we develop a descriptive model, which predicts the best time for inspections and find the optimum values of predefined decision variables such as reliability, preventive maintenance threshold, upper limit for spare part ordering, reorder level of deterioration, and successive imperfect maintenance threshold. The outstanding feature of this research is simultaneous application of bellow concepts for optimization of spare parts inventory and maintenance planning.

1- Using imperfect maintenance
2- Designing an innovative policy for ordering spare parts which is developed based on deterioration characteristic of the system.
3- Employing an integrated approach, which considers maintenance and inventory policies simultaneously, for determining inspection time.
4- Considering multi-location spare part supply chain.
5- Implementing opportunistic, predictive, and preventive maintenance approaches simultaneously.

Simultaneous consideration of all above concepts, distinguishes this research from the other ones none of which have developed an integrated model like this. For future studies, we recommend including these new ideas to develop more comprehensive models.

- ✓ multi-echelon and multi-indenture supply chain
- ✓ the correlation between sub-systems
- ✓ multi-component systems
- ✓ The constraints which are related to limitation of maintenance operations such as manpower, quality of operations and the number of equipment.
- ✓ the uncertainty in the spare parts procurement
- ✓ the competitive price among the local suppliers
- ✓ the environmental parameters and standards
- ✓ the effects of maintenance on the quality of production

Considering all above items helps us develop more realistic models, but these new assumptions make the problems much more complicated. Employing modern solving methods can lead decision makers to find the optimized solution much easier.